\documentclass[11pt,a4paper]{article}

\usepackage[T1]{fontenc}
\usepackage[utf8]{inputenc}
\usepackage[margin=1.05in]{geometry}
\usepackage{amsmath,amssymb,amsthm,mathtools}
\usepackage[hidelinks,pdfusetitle]{hyperref}
\hypersetup{pdftitle={The Dimension-Shift Category and Its Mellin-Gamma Representation},
            pdfauthor={Andreu Ballus Santacana}}

\theoremstyle{plain}
\newtheorem{theorem}{Theorem}[section]
\newtheorem{proposition}[theorem]{Proposition}

\newtheorem{corollary}[theorem]{Corollary}

\theoremstyle{definition}
\newtheorem{definition}[theorem]{Definition}

\theoremstyle{remark}
\newtheorem{remark}[theorem]{Remark}

\newcommand{\R}{\mathbb{R}}

\newcommand{\Z}{\mathbb{Z}}

\newcommand{\Cc}{C_c}
\newcommand{\Dim}{\mathsf{Dim}_+}
\newcommand{\RadMeas}{\mathsf{RadMeas}}
\newcommand{\Mfun}{\mathcal{M}}
\newcommand{\Bobs}{\mathcal{B}}

\title{The Dimension-Shift Category and Its Mellin--Gamma Representation}

\author{Andreu Ballús Santacana\thanks{\texttt{andreu.ballus@uab.cat}.
Departament de Filosofia, Universitat Autònoma de Barcelona; ToposCircuitry, S.L.,
Granollers, Barcelona, Spain.}}

\date{\today}

\begin{document}

\maketitle

\begin{center}
\fbox{\begin{minipage}{0.92\textwidth}
\small
\textbf{Note on v2.} Substantial revision of v1 (\emph{Analytic Uniqueness of Ball Volume Interpolation}, June 2025). The categorical project initiated in v1 has split into two complementary papers: this one (categorical) and the analytic companion \emph{Radial Integration and Ball Volume in Continuous Dimension}~\cite{BallusAnalytic} (arXiv:2605.04351, math.CA). The main result is now a functor classification on a thin category $\Dim$ of dimension shifts; the categorical-dimension reading $\beta(x) = \dim_{\mathrm{cat}}(\mathsf X_x) = x$ is a consequence of the classification, not a hypothesis. The unitary and symplectic statements of v1 are not asserted here.
\end{minipage}}
\end{center}

\bigskip

\begin{abstract}
We define a thin category $\Dim$ of dimension shifts on the positive reals and a category $\RadMeas$ of positive Radon measures on $(0,\infty)$ with positive Radon--Nikodym densities as morphisms. We classify the scaling-covariant functors $\Dim\to\RadMeas$ whose morphism component is given by homogeneous densities: such functors are parameterized by their object-level coefficient function, and the cocycle condition for dimension shifts coincides with the functoriality axiom. Imposing Gaussian normalization on objects then selects a unique functor $\Mfun$, whose values are the Mellin--Gamma measures $d\mu_x(u) = \frac{\pi^{x/2}}{\Gamma(x/2)}u^{x/2-1}\,du$. The morphism part of $\Mfun$ recovers the radial-integration transport $R(x,r) = \pi^r\,\Gamma(x/2)/\Gamma(x/2+r)$. Composing $\Mfun$ with the unit-interval observable $\Bobs(\mu) := \mu((0,1))$ yields a scalar transport $T(x,r)$ that differs from $R$ by the multiplicative coboundary of $\beta(x) = x$; this coboundary is the categorical-dimension functional of the standard object in Deligne's interpolation category $\mathrm{Rep}(O_t)$. As a corollary, the Euclidean ball-volume formula $V(x) = \pi^{x/2}/\Gamma(x/2+1)$ is the value of $\Bobs\circ\Mfun$ at $x$. The main object of the paper is the functor $\Mfun$; the ball-volume formula is one of its scalar evaluations.
\end{abstract}

\section{Introduction}\label{sec:intro}

This paper is about a functor. Specifically, we classify functors from a category of dimension shifts to a category of positive Radon measures on the positive half-line, subject to two axioms: scaling covariance of objects and a homogeneous-density form for morphisms. The classification is the central result. As a scalar consequence, the Euclidean ball-volume formula $V(x) = \pi^{x/2}/\Gamma(x/2+1)$ in continuous dimension $x>0$ emerges as the value of an observable on the unique functor selected by a Gaussian normalization condition.

The categorical organization makes visible a structure that is not visible at the level of individual measures. The radial Mellin--Gamma family
\[
d\mu_x(u) \;=\; \frac{\pi^{x/2}}{\Gamma(x/2)}\,u^{x/2-1}\,du, \qquad x>0,
\]
when treated as the object part of a functor $\Mfun\colon\Dim\to\RadMeas$, has a forced morphism part: for every dimension shift $r\colon x\to x+2r$, the Radon--Nikodym density of $\mu_{x+2r}$ with respect to $\mu_x$ is the homogeneous function
\[
h_{x,r}(u) \;=\; R(x,r)\,u^r, \qquad R(x,r) \;=\; \pi^r\,\frac{\Gamma(x/2)}{\Gamma(x/2+r)}.
\]
Functoriality of $\Mfun$ is the cocycle identity $R(x,r+s) = R(x+2r,s)\,R(x,r)$, which holds by direct cancellation of intermediate Gamma factors. Postcomposition with the unit-interval observable $\Bobs(\mu) := \mu((0,1))$ yields a scalar transport $T(x,r) = V(x+2r)/V(x)$ that differs from $R$ by the multiplicative coboundary of $\beta(x) = x$:
\[
\frac{R(x,r)}{T(x,r)} \;=\; \frac{x+2r}{x} \;=\; \frac{\beta(x+2r)}{\beta(x)}.
\]
The function $\beta$ is interpretable: in Deligne's rigid symmetric monoidal category $\mathrm{Rep}(O_t)$, the standard generating object $\mathsf X_t$ has categorical dimension equal to the parameter, $\dim_{\mathrm{cat}}(\mathsf X_t) = t$. Thus $\beta(x) = \dim_{\mathrm{cat}}(\mathsf X_x)$, and the coboundary distinguishing the radial-integration transport from the ball-volume transport is the categorical-dimension functional of the standard object.

\paragraph{Structure of the paper.}
Section~\ref{sec:dimcat} defines the dimension-shift category $\Dim$. Section~\ref{sec:radmeas} defines the category $\RadMeas$ of positive Radon measures with Radon--Nikodym density morphisms. Section~\ref{sec:functor} defines scaling-covariant functors with homogeneous-density morphism part and proves the functor classification theorem (Theorem~\ref{thm:functor-classification}). Section~\ref{sec:gauss} imposes the Gaussian normalization condition and identifies the unique resulting functor $\Mfun$ (Theorem~\ref{thm:gauss-functor}). Section~\ref{sec:cocycles} interprets the morphism part of $\Mfun$ as the radial-integration cocycle $R$, the postcomposition with the unit-interval observable as the ball-volume cocycle $T$, and the difference $R/T$ as the coboundary of the linear function $\beta(x)=x$. Section~\ref{sec:catdim} interprets $\beta$ as the categorical-dimension functional in $\mathrm{Rep}(O_t)$ (Proposition~\ref{prop:catdim}). Section~\ref{sec:ball} derives the Euclidean ball-volume formula as a corollary. Section~\ref{sec:scope} records the limits of the categorical content of this paper and points to a companion analytic treatment.

\paragraph{Scope.}
The category theory used here is light: thin categories, functor classification, cocycle conditions, multiplicative coboundaries, and a categorical-trace identity from the construction of $\mathrm{Rep}(O_t)$. We use standard terminology for categories, functors, natural transformations, and rigid symmetric monoidal categories; Riehl~\cite{Riehl2016} provides a self-contained treatment. The paper does not use $2$-categorical structure, Kan extensions, enriched or internal categories, or higher-categorical methods. The contribution is conceptual: it shows that the radial Mellin--Gamma measure is canonically a functor, that the ball-volume formula is a scalar evaluation of this functor under a specific observable, and that the gap between the two natural transports is the coboundary of a categorical-dimension functional.

\paragraph{Companion paper.}
The same mathematical content can be developed analytically, with the classification stated for individual positive functionals on $\Cc(\R_{>0})$ and the functorial structure deduced afterwards. That treatment, together with an independent shifted Bohr--Mollerup characterization of $T$, an embedding of the ball-volume observable into a wider system of Mellin observables, and a homogeneous-radial extension to the Bui--Randles polar framework~\cite{BuiRandles2022}, appears in~\cite{BallusAnalytic}. The present paper takes the categorical organization as primary and refers to the companion for these analytic complements.

\section{The dimension-shift category}\label{sec:dimcat}

\begin{definition}\label{def:dimcat}
The \emph{dimension-shift category} $\Dim$ is the category whose objects are the positive real numbers,
\[
\mathrm{Ob}(\Dim) \;=\; \R_{>0},
\]
and whose morphisms are
\[
\Dim(x,y) \;=\; \begin{cases} \{r\} & \text{if } y = x + 2r \text{ for some } r\ge 0,\\ \emptyset & \text{otherwise.}\end{cases}
\]
Composition is addition of shift parameters: if $r\colon x \to x+2r$ and $s\colon x+2r\to x+2r+2s$, then $s\circ r \;=\; r+s\colon x\to x+2(r+s)$. The identity at $x$ is $0\colon x\to x$.
\end{definition}

The category $\Dim$ is thin (at most one morphism between any two objects), has all objects connected by morphisms in a single direction (positive shifts only), and is the delooping of the additive monoid $\R_{\ge 0}$ acted on by translation, restricted to positive starting objects. Negative shifts $r\colon x\to x+2r$ with $x+2r > 0$ can be incorporated by passing to the appropriate localization; we work with the nonnegative direction throughout.

\begin{remark}\label{rmk:dimcat-meaning}
The name reflects the intended interpretation: an object $x$ is to be thought of as a continuous dimension parameter, and a morphism $r\colon x\to x+2r$ is a shift in dimension by $2r$. The factor of $2$ encodes the choice of working in the squared-radius variable $u=r^2$ in the Euclidean setting; in particular, integer dimensions $x = n$ correspond to integer-step shifts $n \to n+2$. The category $\Dim$ does not itself prefer this interpretation — it is simply the delooping of $\R_{\ge 0}$ — but the geometric content of the paper makes it the natural domain for radial-measure functors.
\end{remark}

\section{The category of positive Radon measures}\label{sec:radmeas}

\begin{definition}\label{def:radmeas}
The category $\RadMeas$ has as objects the nonzero, locally finite positive Radon measures on $(0,\infty)$. A morphism $\mu \to \nu$ is a positive measurable function $h\colon(0,\infty)\to[0,\infty)$ such that
\[
d\nu \;=\; h\,d\mu, \qquad\text{i.e.,}\qquad \nu(A) \;=\; \int_A h\,d\mu \quad\text{for every Borel } A\subset(0,\infty).
\]
Composition is pointwise multiplication of densities: if $h\colon\mu\to\nu$ and $g\colon\nu\to\rho$, then $g\circ h\colon\mu\to\rho$ is the function $u\mapsto g(u)h(u)$. The identity at $\mu$ is the constant function $1$. We identify two morphisms $h,h'\colon\mu\to\nu$ when $h = h'$ $\mu$-almost everywhere.
\end{definition}

The category $\RadMeas$ is well defined: composition is associative because pointwise multiplication of nonnegative measurable functions is associative, and the identity laws hold trivially. The almost-everywhere identification ensures that the morphism set is the correct quotient — Radon--Nikodym densities are determined only up to a $\mu$-null set, so equality of densities pointwise would be too strict.

\begin{remark}\label{rmk:radmeas-existence}
A morphism $\mu\to\nu$ exists if and only if $\nu$ is absolutely continuous with respect to $\mu$, in which case the morphism is unique (modulo the almost-everywhere identification) by the Radon--Nikodym theorem. Hence $\RadMeas$ is essentially a partial-order-with-density structure on the class of mutually absolutely continuous measures.
\end{remark}

\section{Scaling-covariant homogeneous functors}\label{sec:functor}

We now define the class of functors $\Dim\to\RadMeas$ that the paper classifies.

\begin{definition}\label{def:functor-class}
A functor $F\colon\Dim\to\RadMeas$ is called \emph{scaling-covariant homogeneous} if both:

(SC) For every $x>0$, the measure $F(x)$ on $(0,\infty)$ satisfies
\[
\int_0^\infty \phi(\lambda u)\,dF(x)(u) \;=\; \lambda^{-x/2}\int_0^\infty \phi(u)\,dF(x)(u)
\]
for every $\phi\in\Cc(\R_{>0})$ and every $\lambda>0$.

(HD) For every $x>0$ and every $r\ge 0$, the morphism $F(r)\colon F(x)\to F(x+2r)$ is given by a Radon--Nikodym density of homogeneous form
\[
h_{x,r}(u) \;=\; A(x,r)\,u^r
\]
for some positive function $A\colon\{(x,r) : x>0,\,r\ge 0\}\to(0,\infty)$.
\end{definition}

Condition (SC) is a homogeneity statement on objects: in the squared-radial variable $u$, the integration operator represented by $F(x)$ has homogeneity degree $x/2$. Condition (HD) is a homogeneity statement on morphisms: dimension shifts act by multiplication by powers of $u$, with a coefficient depending on $(x,r)$.

\begin{theorem}[Functor classification]\label{thm:functor-classification}
Let $F\colon\Dim\to\RadMeas$ be a scaling-covariant homogeneous functor. Then there exists a positive function $c\colon\R_{>0}\to(0,\infty)$ such that
\begin{equation}\label{eq:Fobj}
dF(x)(u) \;=\; c(x)\,u^{x/2-1}\,du \qquad\text{for every } x>0,
\end{equation}
and the morphism part is forced by
\begin{equation}\label{eq:Fmorph}
A(x,r) \;=\; \frac{c(x+2r)}{c(x)} \qquad\text{for every } x>0,\ r\ge 0.
\end{equation}
Conversely, given any positive function $c$, the formulas~\eqref{eq:Fobj} and~\eqref{eq:Fmorph} define a scaling-covariant homogeneous functor $F\colon\Dim\to\RadMeas$.
\end{theorem}

\begin{proof}
\emph{Object part.} Fix $x>0$. By the Riesz--Markov--Kakutani theorem~\cite[Thm.~7.2]{Folland1999}, the positive linear functional on $\Cc(\R_{>0})$ defined by $\phi\mapsto \int\phi\,dF(x)$ is represented uniquely by the Radon measure $F(x)$. Condition (SC) translates, via the dilation $S_\lambda(u)=\lambda u$, into the identity of measures
\[
(S_\lambda)_*F(x) \;=\; \lambda^{-x/2}F(x)
\]
on $\Cc(\R_{>0})$, hence (by Riesz--Markov uniqueness) on all Borel sets of compact closure in $(0,\infty)$, hence on all Borel sets. Equivalently, $F(x)(\lambda A) = \lambda^{x/2}F(x)(A)$.

Pull back through $q\colon\R\to\R_{>0}$, $q(t)=e^t$, to a Radon measure $\nu_x := (q^{-1})_*F(x)$ on $\R$. The multiplicative scaling law on $F(x)$ becomes the additive quasi-invariance $\nu_x(E+s) = e^{(x/2)s}\nu_x(E)$. The measure $\eta_x$ on $\R$ defined by $d\eta_x(t) = e^{-(x/2)t}\,d\nu_x(t)$ is then translation-invariant on $\R$ (a one-line check via the integrated quasi-invariance relation against $\Cc(\R)$ test functions). By the uniqueness theorem for Haar measure on locally compact groups~\cite[Thm.~11.9]{Folland1999}, there exists $c(x)\ge 0$ with $d\eta_x(t) = c(x)\,dt$, hence $d\nu_x(t) = c(x)\,e^{(x/2)t}\,dt$. Pulling back through $u = e^t$ gives
\[
dF(x)(u) \;=\; c(x)\,u^{(x/2)-1}\,du.
\]
Since $F(x)$ is a nonzero object of $\RadMeas$ by Definition~\ref{def:radmeas}, $c(x) > 0$.

\emph{Morphism part.} For a morphism $r\colon x\to x+2r$, the Radon--Nikodym density $h_{x,r}$ of $F(x+2r)$ with respect to $F(x)$ is determined by~\eqref{eq:Fobj}:
\[
h_{x,r}(u) \;=\; \frac{dF(x+2r)}{dF(x)}(u) \;=\; \frac{c(x+2r)\,u^{(x+2r)/2-1}}{c(x)\,u^{x/2-1}} \;=\; \frac{c(x+2r)}{c(x)}\,u^r.
\]
This matches the homogeneous form (HD) with coefficient $A(x,r) = c(x+2r)/c(x)$, proving~\eqref{eq:Fmorph}.

\emph{Functoriality.} For composable shifts $r\colon x\to x+2r$ and $s\colon x+2r\to x+2(r+s)$, functoriality of $F$ requires $F(s)\circ F(r) = F(r+s)$, i.e.,
\[
A(x+2r,s)\,u^s \cdot A(x,r)\,u^r \;=\; A(x,r+s)\,u^{r+s},
\]
which reduces to $A(x+2r,s)\,A(x,r) = A(x,r+s)$. Substituting~\eqref{eq:Fmorph},
\[
\frac{c(x+2r+2s)}{c(x+2r)}\cdot\frac{c(x+2r)}{c(x)} \;=\; \frac{c(x+2r+2s)}{c(x)} \;=\; A(x,r+s),
\]
so functoriality holds automatically and imposes no further restriction on $c$.

\emph{Converse.} Given any positive $c\colon\R_{>0}\to(0,\infty)$, the assignments~\eqref{eq:Fobj} on objects and~\eqref{eq:Fmorph} on morphisms define a functor: identity preservation is $A(x,0) = c(x)/c(x) = 1$, and composition was verified in the previous paragraph.
\end{proof}

\begin{remark}\label{rmk:no-extra-cocycle}
The functor classification has no residual cocycle condition because the cocycle condition for a multiplicative shift function — $A(x,r+s) = A(x+2r,s)\,A(x,r)$ — is automatic for any function of the form $A(x,r) = c(x+2r)/c(x)$. Equivalently: every multiplicative $1$-cocycle for the dimension-shift action on $\R_{>0}$ that is a coboundary is realized as the morphism part of some scaling-covariant homogeneous functor, and the realization is unique once the object part is fixed.
\end{remark}

\section{Gaussian normalization selects the Mellin--Gamma functor}\label{sec:gauss}

The classification of Theorem~\ref{thm:functor-classification} leaves the coefficient function $c$ free. We now impose a single normalization condition that fixes it.

\begin{definition}\label{def:gauss-norm}
A scaling-covariant homogeneous functor $F\colon\Dim\to\RadMeas$ is \emph{Gaussian-normalized} if for every $x>0$, the function $u\mapsto e^{-u}$ is $F(x)$-integrable on $(0,\infty)$ and
\[
\int_0^\infty e^{-u}\,dF(x)(u) \;=\; \pi^{x/2}.
\]
\end{definition}

\begin{theorem}[Gaussian-normalized functor]\label{thm:gauss-functor}
There is a unique scaling-covariant homogeneous Gaussian-normalized functor
\[
\Mfun\colon\Dim\to\RadMeas.
\]
Its values are
\begin{equation}\label{eq:Mobj}
d\Mfun(x)(u) \;=\; \frac{\pi^{x/2}}{\Gamma(x/2)}\,u^{x/2-1}\,du, \qquad x>0,
\end{equation}
and its morphism part is given by
\begin{equation}\label{eq:Mmorph}
\Mfun(r)(u) \;=\; R(x,r)\,u^r, \qquad R(x,r) \;=\; \pi^r\,\frac{\Gamma(x/2)}{\Gamma(x/2+r)}.
\end{equation}
\end{theorem}

\begin{proof}
By Theorem~\ref{thm:functor-classification}, $dF(x)(u) = c(x)\,u^{x/2-1}\,du$ for some positive $c$. The Gaussian normalization condition gives
\[
\pi^{x/2} \;=\; c(x)\int_0^\infty e^{-u}\,u^{x/2-1}\,du \;=\; c(x)\,\Gamma(x/2),
\]
using the Euler integral representation of the Gamma function~\cite[Ch.~1]{Artin1964},~\cite[\S 4]{FlajoletGourdonDumas1995}, so $c(x) = \pi^{x/2}/\Gamma(x/2)$. The object formula~\eqref{eq:Mobj} follows. The morphism formula~\eqref{eq:Mmorph} follows from~\eqref{eq:Fmorph}:
\[
A(x,r) \;=\; \frac{\pi^{(x+2r)/2}/\Gamma((x+2r)/2)}{\pi^{x/2}/\Gamma(x/2)} \;=\; \pi^r\,\frac{\Gamma(x/2)}{\Gamma(x/2+r)} \;=\; R(x,r). \qedhere
\]
\end{proof}

\begin{remark}[On the choice of normalizing probe]\label{rmk:probe}
The Gaussian probe $u\mapsto e^{-u}$ in Definition~\ref{def:gauss-norm} is one valid choice among many. Any positive nonzero test function $\psi$ for which $\int_0^\infty\psi(u)\,u^{x/2-1}\,du$ converges absolutely and is positive at all relevant $x$ would suffice to fix $c$, with a different prescribed value $\int\psi\,dF(x) = N(x)>0$ on the right-hand side simply yielding a different constant. The Gaussian choice is privileged not by structural uniqueness within the present framework but by its agreement with the classical Euclidean Gaussian integral $\int_{\R^n}e^{-|y|^2}\,dy = \pi^{n/2}$ at integer $x=n$, which makes the resulting functor $\Mfun$ the canonical extension of integer-dimensional Euclidean radial integration. The classification of $\Mfun(x)$ as a Mellin--Gamma density does not depend on this choice; only the specific constant $\pi^{x/2}/\Gamma(x/2)$ does.
\end{remark}

\section{The radial-integration cocycle, the ball-volume cocycle, and their coboundary}\label{sec:cocycles}

The functor $\Mfun$ produces, automatically, the two transports that play the central role in the analytic treatment.

\subsection*{The radial-integration transport as functoriality}

The morphism component $R$ of $\Mfun$, defined by~\eqref{eq:Mmorph}, satisfies the cocycle identity
\[
R(x,r+s) \;=\; R(x+2r,s)\,R(x,r) \qquad\text{for all admissible } x,r,s,
\]
because functoriality of $\Mfun$ requires it (Theorem~\ref{thm:functor-classification}, functoriality verification). Equivalently, viewing $R$ as a function on the morphisms of $\Dim$ valued in positive scalars (via the multiplicative coefficient in the homogeneous density), $R$ is a multiplicative $1$-cocycle on $\Dim$. We refer to $R$ as the \emph{radial-integration cocycle}.

\subsection*{The ball-volume cocycle as a scalar observable}

Define the \emph{unit-interval observable} $\Bobs\colon\RadMeas\to\R_{>0}$ by
\[
\Bobs(\mu) \;:=\; \mu\bigl((0,1)\bigr).
\]
This is not a functor of $\RadMeas$ in any natural sense (the target $\R_{>0}$ is not the morphism set of a category in this paper, and $\Bobs$ does not respect Radon--Nikodym morphisms in a covariant way), but it is a positive real-valued function on objects, and that is all we need.

Define
\[
V(x) \;:=\; \Bobs(\Mfun(x)) \;=\; \Mfun(x)\bigl((0,1)\bigr),
\]
and the associated scalar transport
\[
T(x,r) \;:=\; \frac{V(x+2r)}{V(x)}.
\]
By direct computation from~\eqref{eq:Mobj},
\[
V(x) \;=\; \frac{\pi^{x/2}}{\Gamma(x/2)}\int_0^1 u^{x/2-1}\,du \;=\; \frac{\pi^{x/2}}{\Gamma(x/2)}\cdot\frac{2}{x} \;=\; \frac{\pi^{x/2}}{\Gamma(x/2+1)},
\]
and
\[
T(x,r) \;=\; \pi^r\,\frac{\Gamma(x/2+1)}{\Gamma(x/2+r+1)}.
\]
The function $T$ also satisfies the cocycle identity by direct cancellation of intermediate Gamma factors. We refer to $T$ as the \emph{ball-volume cocycle}.

\subsection*{The coboundary distinguishing $R$ and $T$}

\begin{theorem}[Coboundary identification]\label{thm:coboundary}
The cocycles $R$ and $T$ differ by the multiplicative coboundary of the function $\beta\colon\R_{>0}\to(0,\infty)$, $\beta(x) = x$:
\[
\frac{R(x,r)}{T(x,r)} \;=\; \frac{x+2r}{x} \;=\; \frac{\beta(x+2r)}{\beta(x)} \qquad\text{for all admissible } x,r.
\]
Equivalently, $T(x,r) = \frac{\beta(x)}{\beta(x+2r)}\,R(x,r)$.
\end{theorem}

\begin{proof}
Direct from the closed forms in Theorem~\ref{thm:gauss-functor} and the previous subsection, using $\Gamma(x/2+1) = (x/2)\Gamma(x/2)$ and $\Gamma(x/2+r+1) = (x/2+r)\Gamma(x/2+r)$:
\begin{align*}
\frac{R(x,r)}{T(x,r)}
&= \frac{\pi^r\,\Gamma(x/2)/\Gamma(x/2+r)}{\pi^r\,\Gamma(x/2+1)/\Gamma(x/2+r+1)}
= \frac{\Gamma(x/2)\,\Gamma(x/2+r+1)}{\Gamma(x/2+r)\,\Gamma(x/2+1)} \\
&= \frac{x/2+r}{x/2} = \frac{x+2r}{x}. \qedhere
\end{align*}
\end{proof}

The pair $(R, T)$ therefore determines a single cohomology class in the multiplicative cocycle theory of $\Dim$, and the coboundary $\beta(x) = x$ is the algebraic record of the gap between functorial transport (which $R$ realizes strictly) and observable-induced transport (which $T$ realizes after applying the unit-interval observable). The function $V(x)$ on the objects of $\Dim$ is therefore not preserved strictly under dimension shift: it is preserved up to the coboundary defect of $\beta$.

\section{The categorical-dimension interpretation of the coboundary}\label{sec:catdim}

The function $\beta(x) = x$ that controls the coboundary in Theorem~\ref{thm:coboundary} admits a natural reading in Deligne's interpolation categories. We make this precise as follows.

\begin{proposition}\label{prop:catdim}
Let $\mathrm{Rep}(O_t)$ denote Deligne's rigid symmetric monoidal interpolation category for the orthogonal series; the framework is recalled in~\cite[\S\S 9--10]{Deligne2007}, and we use the Brauer-category presentation of~\cite[\S\S 2.1--2.2]{ComesHeidersdorf2017}, in which $\mathrm{Rep}(O_t)$ is the Karoubi envelope of the additive envelope of the Brauer category $\mathcal B(t)$. Let $\mathsf X_t$ denote its standard generating object. Then $\dim_{\mathrm{cat}}(\mathsf X_t) = t$, and consequently
\[
\frac{R(x,r)}{T(x,r)} \;=\; \frac{\dim_{\mathrm{cat}}(\mathsf X_{x+2r})}{\dim_{\mathrm{cat}}(\mathsf X_x)} \qquad\text{for all admissible } x,r.
\]
That is, the coboundary distinguishing the radial-integration cocycle from the ball-volume cocycle is the categorical-dimension functional of the standard object in $\mathrm{Rep}(O_t)$.
\end{proposition}

\begin{proof}
The categorical dimension of an object in a rigid symmetric monoidal category is the categorical trace of its identity morphism: the scalar in $\mathrm{End}(\mathbf 1)$ obtained by composing coevaluation $\mathbf 1\to X\otimes X^*$ with evaluation $X^*\otimes X\to\mathbf 1$ via the symmetry. In the Brauer category $\mathcal B(t)$, this trace on the generating object is the closed loop, which the defining relations of $\mathcal B(t)$ assign the scalar $t$~\cite[\S 2.1]{ComesHeidersdorf2017}; compare~\cite[\S 9.2]{Deligne2007}. Categorical dimension is preserved under additive envelope and Karoubi envelope, so $\dim_{\mathrm{cat}}(\mathsf X_t) = t$ in $\mathrm{Rep}(O_t)$. Combining with Theorem~\ref{thm:coboundary} gives the stated identity.
\end{proof}

\begin{remark}[Scope of the categorical interpretation]\label{rmk:catdim-scope}
Proposition~\ref{prop:catdim} does not derive any of the analytic content of the paper from $\mathrm{Rep}(O_t)$. The functor $\Mfun$ is constructed in Sections~\ref{sec:functor}--\ref{sec:gauss} from the axioms (SC), (HD), and Gaussian normalization, with no reference to interpolation categories. Proposition~\ref{prop:catdim} adds an external interpretation of the dimension parameter: an object $x$ of $\Dim$ can be read as the categorical dimension of a canonical object in a rigid symmetric monoidal category that exists for independent reasons. In this reading, the codomain of $\beta\colon\mathrm{Ob}(\Dim)\to\R_{>0}$ is the categorical-dimension functional, and the coboundary $\beta(x+2r)/\beta(x)$ is its multiplicative shift under a dimension shift of $2r$. The interpretation is structural rather than load-bearing.
\end{remark}

\section{The Euclidean ball-volume formula as a corollary}\label{sec:ball}

\begin{corollary}[Euclidean ball volume]\label{cor:ball-volume}
The composite $\Bobs\circ\Mfun\colon\mathrm{Ob}(\Dim)\to\R_{>0}$ takes the value
\[
V(x) \;=\; \frac{\pi^{x/2}}{\Gamma(x/2+1)}, \qquad x>0.
\]
At integer $x=n\in\Z_{\ge 1}$, $V(n)$ equals the classical Lebesgue volume of the unit ball in $\R^n$.
\end{corollary}

\begin{proof}
The closed form was computed in Section~\ref{sec:cocycles}. For integer $x=n$, the change of variable $u=r^2$ relates the unit interval $u\in(0,1)$ to the unit ball $r\in(0,1)$, and~\eqref{eq:Mobj} at $x=n$ recovers the polar-coordinate factor $\pi^{n/2}/\Gamma(n/2)$ for the Lebesgue radial measure on $\R^n$ in $u$-coordinates~\cite[\S 3]{BallusAnalytic}. The identity $V(n) = \pi^{n/2}/\Gamma(n/2+1)$ is then the standard Euclidean ball-volume formula.
\end{proof}

This is the functorial origin of the formula. The Euclidean ball-volume formula is the value of a single canonical observable ($\Bobs$) on a single canonical functor ($\Mfun$), with no recourse to interpolation, holomorphic continuation, or analytic special-function arguments.

\section{Scope, limits, and companion paper}\label{sec:scope}

The categorical content of this paper is light by design. The main objects ($\Dim$ and $\RadMeas$) are 1-categories of an elementary kind: $\Dim$ is the delooping of an additive monoid restricted to positive starting objects, and $\RadMeas$ is essentially a partial-order with morphism data given by Radon--Nikodym densities. The main constructions (scaling covariance, homogeneous-density morphism part, Gaussian normalization) are functorial conditions that do not invoke higher-categorical structure, enrichment, monoidal coherence, adjunctions, or Kan extensions. The categorical-dimension functional appearing in Proposition~\ref{prop:catdim} is the only nontrivial categorical input from outside, and it enters only as an interpretation of the coboundary, not as part of the construction of $\Mfun$.

What the paper does contribute, categorically, is the recognition that the Mellin--Gamma family of measures is not just a parameterized family but a specific functor with a specific morphism part — and that this morphism part is what the analytic transport $R(x,r)$ realizes. The cocycle $T$ is then visible as a non-functorial scalar quantity obtained by postcomposing with an observable that does not respect the morphism structure, and the gap between $R$ and $T$ is the coboundary of the dimension function. This is a categorical organization of structure that already exists analytically; the contribution is the organization, not new mathematical content.

For an analytic treatment in which the classification is stated for individual positive functionals on $\Cc(\R_{>0})$ rather than for functors on $\Dim$, together with an independent shifted Bohr--Mollerup characterization of $T$, an embedding of the ball-volume observable into a wider family of Mellin observables (sublevel volumes and Gaussian moment observables), and a homogeneous-radial extension to the Bui--Randles polar framework, see~\cite{BallusAnalytic}.

\section*{AI use disclosure}
In preparing this manuscript the author made use of several large language models as technical-editorial interlocutors, in particular for prose tightening, structural feedback on the categorical organization, and verification passes on argument flow. All mathematical content---definitions, theorem statements, proofs, and the choice of structural framework---is the author's. The author has independently verified every claim in the paper and accepts full responsibility for its content.

\end{document}